\documentclass[a4paper,11pt]{amsart}
\usepackage[
  top=3cm,
  bottom=3cm,
  left=2cm,
  right=2cm,
  headheight=14pt,
  headsep=1cm,
  footskip=1cm
]{geometry}
\usepackage{amsmath, epsfig, verbatim}
\usepackage{amsfonts}
\usepackage{amsthm}
\usepackage{mathtools}
\usepackage{amsmath}
\usepackage{amssymb} 
\usepackage{graphics,graphicx}
\usepackage{epstopdf}
\usepackage{listings}
\usepackage{xcolor}
\usepackage{tikz}
\usepackage{subcaption}
\usepackage{pgfplots}
\usepackage{caption}
\usepackage{booktabs}
\usepackage{tabularx}
\usepackage{hyperref}   
\usepackage{doi}
\usepackage{enumitem}
\usepackage[numbers,sort&compress]{natbib}

\newtheorem{theorem}{Theorem}[section]
\newtheorem{lemma}[theorem]{Lemma}
\newtheorem{prop}[theorem]{Proposition}

\theoremstyle{definition}

\theoremstyle{remark}

\begin{document}
\title{Mutual visibility in $(d, 5)$-graphs with excess}
\author{Tonny K B}
\address{Tonny K B, Department of Mathematics, College of Engineering Trivandrum, Thiruvananthapuram, Kerala, India, 695016.}
\email{tonnykbd@cet.ac.in}

\begin{abstract}
The author investigates mutual visibility in regular graphs of girth
at least $5$, with particular emphasis on graphs with excess. A general upper bound on the mutual-visibility number is established
in terms of the degree and order. More significantly, an algebraic
inequality is derived that restricts the mutual-visibility number of
$(d,5,+\delta)$-graphs. As applications, the mutual-visibility number
of the Robertson graph is shown to be $9$, and upper bounds are
obtained for the $(5,5)$-cages and the $(6,5)$-cage.
\end{abstract}
\maketitle

\section{Introduction}
Let $G$ be a connected graph and let $S\subseteq V(G)$.
Two vertices $u,v\in S$ are said to be \emph{$S$-visible} \cite{MV_3}
if there exists a shortest $(u,v)$-path whose internal
vertices lie outside $S$. The set $S$ is a
\emph{mutual-visibility set} \cite{Stefano} if every two distinct
vertices of $S$ are $S$-visible. The maximum
cardinality of a mutual-visibility set of $G$
is called the \emph{mutual-visibility number} of $G$,
denoted by $\mu(G)$.

A $d$-regular graph of girth $g$ is called a \emph{$(d,g)$-graph}.
Sachs~\cite{Sach_1963} proved that $(d,g)$-graphs exist for every
$d\geq 3$ and $g\geq 5$. The degree-girth problem asks for the
minimum order of a $(d,g)$-graph. A graph attaining this minimum is
called a \emph{$(d,g)$-cage}. The order of a $(d,2r+1)$-graph, where $r\geq 2$, is bounded below by the Moore bound
$M(d,2r+1)=1+d\sum_{i=0}^{r-1}(d-1)^i$. Let $G$ be a $(d,2r+1)$-graph of order $n$, and define its \emph{excess} by $\delta=n-M(d,2r+1)\geq 0$. A $(d,2r+1)$-graph with excess $\delta$ is called a $(d,2r+1,+\delta)$-graph. The structural restrictions imposed by the girth
and excess naturally motivate the study of mutual visibility in
this class of graphs.
\section{Notation and Preliminaries}
Throughout this article, $G$ denotes a simple, connected
graph with vertex set $V(G)$ and edge set $E(G)$. The maximum degree of a graph $G$ is denoted by $\Delta(G)$. A sequence of
vertices $(u_0,u_1,\ldots,u_n)$ is called a $(u_0,u_n)$-path in
$G$ if $u_{i-1}u_i\in E(G)$ for every $1\leq i\leq n$. For a vertex $u\in V(G)$, $N_G(u)$ denotes the open neighbourhood of
$u$, and, for an integer $r\geq 0$, $N_G^r[u]$ denotes the set of
vertices at distance at most $r$ from $u$ in $G$.
For $S\subseteq V(G)$, $G[S]$ denotes the subgraph of $G$ induced by $S$. If
$S,T\subseteq V(G)$, then $e_G(S)$ denotes the number of edges of $G[S]$, while $e_G(S,T)$ denotes the number of edges of $G$ with one
end in $S$ and the other in $T$. A function $f$ is \emph{convex} in an interval $[a,b]$ if
$f((1-\lambda)x+\lambda y)\leq(1-\lambda)f(x)+\lambda f(y)$ for all $x,y\in[a,b]$ and $\lambda\in[0,1]$, which in particular
implies that $f(x)\leq\max\{f(a),f(b)\}$ for every $x\in[a,b]$. An \emph{independent set} of a graph $G$ is a set of pairwise nonadjacent vertices. The maximum cardinality of an independent set of $G$ is called the \emph{independence number}
of $G$, denoted by $\alpha(G)$.
\section{Mutual-visibility in $(d,2r+1,+\delta)$-graph}
For a $(d,2r+1,+\delta)$-graph $G$, define the
\emph{excess graph} $\mathcal E_G$ on $V(G)$ by $
uv\in E(\mathcal E_G)$ if and only if $
d_G(u,v)>r$.
\begin{prop}\label{P17.prop1}
Let $G$ be a $(d,2r+1,+\delta)$-graph. Then the excess graph
$\mathcal E_G$ is $\delta$-regular.
\end{prop}

\begin{proof}
Let $u\in V(G)$. Since $G$ has girth $2r+1$, no two distinct paths
of length at most $r$ starting at $u$ can have the same endpoint. Hence
$$|N_G^r[u]|
=1+d\sum_{i=0}^{r-1}(d-1)^i
=M(d,2r+1).$$
Since $|V(G)|=M(d,2r+1)+\delta$, exactly $\delta$ vertices are at
distance greater than $r$ from $u$. By the definition of the excess
graph, these are precisely the neighbours of $u$ in $\mathcal E_G$.
Thus $d_{\mathcal E_G}(u)=\delta$. Since $u$ was arbitrary,
$\mathcal E_G$ is $\delta$-regular.
\end{proof}

The girth condition alone imposes a strong restriction on the
subgraph induced by a mutual-visibility set.
\begin{lemma}\label{P17.lem1}
Let $G$ be a graph of girth at least $5$, and let $S$ be a
mutual-visibility set of $G$. Then $\Delta(G[S])\leq 1$.
Consequently, $2e_G(S)\leq |S|$.
\end{lemma}

\begin{proof}
Suppose that a vertex $x\in S$ has two distinct neighbours
$u,v\in S$. Since $G$ is triangle-free, $u$ and $v$ are non-adjacent,
and hence $(u,x,v)$ is a $(u,v)$-geodesic of length $2$. Moreover,
since $G$ contains no $4$-cycle, $x$ is the unique common neighbour
of $u$ and $v$. Thus $(u,x,v)$ is the unique $(u,v)$-geodesic.
Since its internal vertex $x$ belongs to $S$, this contradicts the
mutual visibility of $S$. Therefore $\Delta(G[S])\leq 1$. By the handshaking lemma,
$2e_G(S)=\sum_{v\in S}d_{G[S]}(v)\leq |S|$, since
$d_{G[S]}(v)\in\{0,1\}$ for every $v\in S$.
\end{proof}
\begin{lemma}[{\cite[Lem.~7]{TonShi_dif2026}}]\label{P6.lem7}
Let $G$ be a $d$-regular graph, let $S\subseteq V(G)$, and put $T=V(G)\setminus S$. Then $e_G(S,T)=d|S|-2e_G(S)$.
\end{lemma}

\begin{theorem}\label{P17.th1}
Let $G$ be a $d$-regular graph of order $n$ and girth at least $5$.
Then
\[
\mu(G)\leq \left\lfloor\frac{dn}{2d-1}\right\rfloor.
\]
\end{theorem}

\begin{proof}
Let $S$ be a mutual-visibility set of $G$, and set
$T=V(G)\setminus S$. By Lemma~\ref{P6.lem7},
$e_G(S,T)=d|S|-2e_G(S)$. Since $2e_G(S)\leq |S|$ by
Lemma~\ref{P17.lem1}, it follows that
$e_G(S,T)\geq (d-1)|S|$. On the other hand, since $G$ is
$d$-regular, $e_G(S,T)\leq d|T|=d(n-|S|)$. Therefore
$(d-1)|S|\leq d(n-|S|)$, and hence
$(2d-1)|S|\leq dn$. Thus $|S|\leq \lfloor dn/(2d-1)\rfloor$. Since $S$ was arbitrary,
the result follows.
\end{proof}
\subsection{Graphs of girth five}
\begin{lemma}\label{P17.lem2}
Let $G$ be a $(d,5,+\delta)$-graph and $S$ a mutual-visibility set.
Set $T=V(G)\setminus S$. Then
\[
\binom{|S|}{2}-e_G(S)-e_{\mathcal E_G}(S)
=
\sum_{w\in T}\binom{|N_G(w)\cap S|}{2}.
\]
\end{lemma}

\begin{proof}
Since $G$ has girth $5$, any two vertices at distance $2$ have a
unique common neighbour. The unordered pairs of vertices in $S$
can be partitioned according as their distance in $G$ is $1$, $2$,
or greater than $2$. The numbers of pairs in the first and third
classes are $e_G(S)$ and $e_{\mathcal E_G}(S)$, respectively.
Consequently, the number of pairs at distance $2$ is
\[
\binom{|S|}{2}-e_G(S)-e_{\mathcal E_G}(S).
\]

Let $\{u,v\}\subseteq S$ with $d_G(u,v)=2$, and let $w$ be their
unique common neighbour. Since $S$ is a mutual-visibility set, the
unique $(u,v)$-geodesic cannot have an internal vertex in $S$.
Hence $w\in T$. For each $w\in T$, every unordered pair of vertices in
$N_G(w)\cap S$ is at distance $2$ in $G$. Furthermore, no such pair
is counted for two distinct vertices of $T$, since two distinct
common neighbours would yield a $4$-cycle. Thus, summing the number of unordered pairs in $N_G(w)\cap S$ over
all $w\in T$, the result follows.
\end{proof}

The identity in Lemma~\ref{P17.lem2} can also be expressed in terms
of the complement of $S$. Let $T=V(G)\setminus S$. Since $G$ is
$d$-regular, Lemma~\ref{P6.lem7} gives
$d|S|=2e_G(S)+e_G(S,T)$ and
$d|T|=2e_G(T)+e_G(S,T)$. Subtracting these equalities and using $|T|=n-|S|$, we obtain
\begin{equation}\label{P17.eq1}
e_G(S)=e_G(T)+\frac{d(2|S|-n)}{2}.    
\end{equation}

Similarly, since $\mathcal E_G$ is $\delta$-regular by
Proposition~\ref{P17.prop1},
\begin{equation}\label{P17.eq2}
e_{\mathcal E_G}(S)
=e_{\mathcal E_G}(T)+\frac{\delta(2|S|-n)}{2}.
\end{equation}
Moreover, $|N_G(w)\cap S|=d-d_{G[T]}(w)$ for every $w\in T$.
Substituting these identities into Lemma~\ref{P17.lem2} yields
\begin{equation}\label{P17.eq3}
 \binom{|S|}{2}
-e_G(T)-e_{\mathcal E_G}(T)
-\frac{(d+\delta)(2|S|-n)}{2}
=
\sum_{w\in T}
\binom{d-d_{G[T]}(w)}{2}.   
\end{equation}

\begin{theorem}\label{P17.th2}
Let $G$ be a $(d,5,+\delta)$-graph with $d\geq 3$, and let $S$ be a
mutual-visibility set of $G$ with $|S|=s$. Then
$s^2-(d+\delta+3)s-(d-3)(d^2+1+\delta)\leq 0$.
\end{theorem}
\begin{proof}
Let $n=|V(G)|=d^2+1+\delta$ and $T=V(G)\setminus S$.
By Lemma~\ref{P17.lem1}, $S\neq V(G)$, since $d\geq 3$.
Hence $|T|=n-s>0$. By Lemma~\ref{P6.lem7},
$e_G(S,T)=ds-2e_G(S)$. Moreover,
\begin{equation}\label{P17.eq20}
\sum_{w\in T}|N_G(w)\cap S|
=e_G(S,T)
=ds-2e_G(S).
\end{equation}

By Lemma~\ref{P17.lem2},
\[
\binom{s}{2}-e_G(S)-e_{\mathcal E_G}(S)
=
\sum_{w\in T}\binom{|N_G(w)\cap S|}{2}.
\]
Since $e_{\mathcal E_G}(S)\geq 0$, we obtain
\begin{equation}\label{P17.eq21}
\binom{s}{2}-e_G(S)
\geq
\sum_{w\in T}\binom{|N_G(w)\cap S|}{2}.
\end{equation}

Using $\binom{x}{2}=(x^2-x)/2$, we have
\[
\sum_{w\in T}\binom{|N_G(w)\cap S|}{2}
=
\frac{1}{2}
\left(
\sum_{w\in T}|N_G(w)\cap S|^2
-
\sum_{w\in T}|N_G(w)\cap S|
\right).
\]
By the Cauchy-Schwarz inequality and \eqref{P17.eq20},
\[
\sum_{w\in T}|N_G(w)\cap S|^2
\geq
\frac{\left(\sum_{w\in T}|N_G(w)\cap S|\right)^2}{|T|}
=
\frac{(ds-2e_G(S))^2}{n-s}.
\]
Consequently,
\begin{equation}\label{P17.eq22}
\sum_{w\in T}\binom{|N_G(w)\cap S|}{2}
\geq
\frac{1}{2}
\left(
\frac{(ds-2e_G(S))^2}{n-s}
-(ds-2e_G(S))
\right).
\end{equation}
Combining \eqref{P17.eq21} and \eqref{P17.eq22} gives
\begin{equation}\label{P17.eq23}
\binom{s}{2}-e_G(S)
\geq
\frac{1}{2}
\left(
\frac{(ds-2e_G(S))^2}{n-s}
-(ds-2e_G(S))
\right).
\end{equation}

\noindent\textbf{Case 1.} $0\leq s\leq n/d$.
Let $f(x)=x^2-(d+\delta+3)x-(d-3)n$. Since $f''(x)=2>0$,
the function $f$ is convex.
Moreover, $f(0)=-(d-3)n\leq 0$. Since
$n=d^2+1+\delta>2d$, we also have
$f\left(\frac{n}{d}\right)
=
\frac{n(d-1)(2d-n)}{d^2}<0$.
Thus, by convexity, $f(s)\leq 0$ for every
$0\leq s\leq n/d$, and the desired inequality follows.

\noindent\textbf{Case 2.}  $s\geq n/d$. For fixed $s$, define
\[
\Phi(t)=
\binom{s}{2}-t
-\frac{1}{2}
\left(
\frac{(ds-2t)^2}{n-s}-(ds-2t)
\right).
\]
By \eqref{P17.eq23}, $\Phi(e_G(S))\geq 0$. Moreover,
\[
\Phi'(t)
=
\frac{2\bigl(ds-2t-(n-s)\bigr)}{n-s}.
\]
For $0\leq t\leq s/2$, we have $ds-2t\geq(d-1)s$.
Since $s\geq n/d$, it follows that $(d-1)s\geq n-s$.
Hence $\Phi'(t)\geq 0$ on $[0,s/2]$, and therefore $\Phi$ is
increasing on this interval. By Lemma~\ref{P17.lem1},
$e_G(S)\leq s/2$, and thus
$0\leq \Phi(e_G(S))
\leq \Phi\left(s/2\right)$.
It follows that
\[
\frac{s(s-2)}{2}
-
\frac{1}{2}
\left(
\frac{(d-1)^2s^2}{n-s}
-(d-1)s
\right)
\geq 0.
\]
Since $s>0$ and $n-s>0$, multiplying by $2(n-s)$ and dividing by
$s$ yields $s^2-(d+\delta+3)s-(d-3)(d^2+1+\delta)\leq 0$, as required.
\end{proof}

\begin{theorem}\label{P17.th3}
Let $G$ be a $(4,5,+2)$-graph. Then $
\mu(G) \leq 9$.
\end{theorem}

\begin{proof}
By Theorem~\ref{P17.th2}, every mutual-visibility set $S$ of $G$
satisfies $|S|^2-9|S|-19\leq 0$, and hence $|S|\leq 10$.
Suppose that $|S|=10$, and set $T=V(G)\setminus S$. Since
$|V(G)|=19$, we have $|T|=9$. By \eqref{P17.eq1},
$e_G(S)=e_G(T)+2$. Moreover, Lemma~\ref{P17.lem1} gives
$\Delta(G[S])\leq 1$, and hence $e_G(S)\leq 5$. Consequently,
$e_G(T)\leq 3$. By \eqref{P17.eq3},
\[
42-e_G(T)-e_{\mathcal E_G}(T)
=
\sum_{w\in T}
\binom{4-d_{G[T]}(w)}{2}.
\]
Using $\binom{4-x}{2}=6-3x+\binom{x}{2}$
and
$\sum_{w\in T}d_{G[T]}(w)=2e_G(T)$,
we obtain
\[
e_{\mathcal E_G}(T)
=
5e_G(T)-12-
\sum_{w\in T}\binom{d_{G[T]}(w)}{2}.
\]
Since both $e_{\mathcal E_G}(T)$ and the sum on the right-hand side
are nonnegative, we obtain $5e_G(T)\geq 12$. As $e_G(T)$ is an
integer and $e_G(T)\leq 3$, it follows that $e_G(T)=3$ and
$e_G(S)=5$. Since $|S|=10$ and $\Delta(G[S])\leq 1$, we conclude
that $G[S]\cong 5K_2$. Moreover, $G[T]$ is triangle-free, since
$G$ has girth $5$. Together with $|T|=9$ and $e_G(T)=3$, this
implies that $G[T]$ is isomorphic to one of the following four
graphs: $3K_2\cup 3K_1$, $P_3\cup K_2\cup 4K_1$,
$P_4\cup 5K_1$, or $K_{1,3}\cup 5K_1$.

Let $w\in T$. The vertex $w$ has
$4-d_{G[T]}(w)$ neighbours in $S$.
Since $G[S]\cong5K_2$ and $G$ is four regular, every vertex of $S$
has exactly three neighbours in $T$.
Thus the number of vertices of $T$
at distance $2$ from $w$ through $S$ is
$2(4-d_{G[T]}(w))$.
Similarly, the number of vertices of $T$
at distance $2$ from $w$ through $T$ is
\[
\sum_{z\in N_{G[T]}(w)}
\bigl(d_{G[T]}(z)-1\bigr).
\]
Since $G$ has girth $5$, these vertices
are all distinct. As $|T|=9$, subtracting
the numbers of vertices at distances
$1$ and $2$ from the eight vertices of
$T\setminus\{w\}$ gives
\begin{equation}\label{P17.eq10}
\begin{aligned}
d_{\mathcal E_G[T]}(w)
&=8-d_{G[T]}(w)
-2(4-d_{G[T]}(w))
-\sum_{z\in N_{G[T]}(w)}
\bigl(d_{G[T]}(z)-1\bigr)\\
&=2d_{G[T]}(w)
-\sum_{z\in N_{G[T]}(w)}d_{G[T]}(z).
\end{aligned}
\end{equation}
If $G[T]\cong K_{1,3}\cup 5K_1$, then, for any leaf $w$ of the
$K_{1,3}$, \eqref{P17.eq10} gives
$d_{\mathcal E_G[T]}(w)=2-3=-1$, a contradiction. If
$G[T]\cong P_4\cup 5K_1$, let $(v_1,v_2,v_3,v_4)$ denote the path
$P_4$. By \eqref{P17.eq10},
$d_{\mathcal E_G[T]}(v_2)=d_{\mathcal E_G[T]}(v_3)=1$, while every
other vertex of $T$ has degree $0$ in $\mathcal E_G[T]$. Hence
$v_2v_3\in E(\mathcal E_G)$, contradicting
$v_2v_3\in E(G)$.

Since $G[S]\cong 5K_2$, we have
$\alpha(G[S])=5$. Let $uv\in E(G[T])$. Since $G$ has girth $5$,
the vertices of $(N_G(u)\cup N_G(v))\cap S$
are pairwise nonadjacent; otherwise, $G$
would contain a triangle or a $4$-cycle.
Moreover, $u$ and $v$ have no common
neighbour in $S$. Since $G$ is $4$-regular, the vertices $u$ and $v$
have $4-d_{G[T]}(u)$ and $4-d_{G[T]}(v)$
neighbours in $S$, respectively. As these
neighbourhoods are disjoint, their union
contains exactly
$(4-d_{G[T]}(u))+(4-d_{G[T]}(v))$ vertices.
Since this union is an independent set in
$G[S]$, its cardinality cannot exceed
$\alpha(G[S])$. Therefore, $(4-d_{G[T]}(u))+(4-d_{G[T]}(v)) \leq 5$.
Consequently, $d_{G[T]}(u)+d_{G[T]}(v)\geq 3$. In particular, $G[T]$ cannot contain a component
isomorphic to $K_2$. This excludes
the cases $G[T]\cong 3K_2\cup 3K_1$ and
$G[T]\cong P_3\cup K_2\cup 4K_1$. Therefore no mutual-visibility set of
size $10$ exists, and hence $\mu(G) \leq 9$.
\end{proof}
The cages of girth $5$ are known for degrees $3\leq d\leq 7$.
Among these, the Petersen, Robertson, and Hoffman-Singleton graphs
are the unique cages of degrees $3$, $4$, and $7$, respectively
\cite{hoffman_moore,Rob_1964}.
There are four nonisomorphic $(5,5)$-cages \cite{Weg_1973},
whereas the $(6,5)$-cage is unique \cite{Keef_Wong_1979, Wong_1979}.
Table~\ref{P17.tab1} summarises the known
mutual-visibility numbers and the upper bounds
obtained from Theorem~\ref{P17.th2} for these cages. The bound in Theorem~\ref{P17.th2} is sharp
for the Petersen and Hoffman-Singleton graphs,
as shown in Table~\ref{P17.tab1}.

\begin{figure}[h]
\centering
\begin{tikzpicture}[ scale= 0.6,
    x=0.026cm,
    y=-0.026cm,
    line width=0.7pt,
    every node/.style={
        circle,
        draw=black,
        fill=white,
        inner sep=0pt,
        minimum size=3.5mm,
        font=\tiny
    }
]

\coordinate (v0) at (124,29);
\coordinate (v1) at (245,29);
\coordinate (v2) at (330,114);
\coordinate (v3) at (330,235);
\coordinate (v4) at (245,320);
\coordinate (v5) at (124,320);
\coordinate (v6) at (38,235);
\coordinate (v7) at (38,114);

\coordinate (v8)  at (107,97);
\coordinate (v9)  at (184,65);
\coordinate (v10) at (262,97);
\coordinate (v11) at (295,175);
\coordinate (v12) at (262,253);
\coordinate (v13) at (184,286);
\coordinate (v14) at (107,253);
\coordinate (v15) at (74,175);

\coordinate (v16) at (171,151);
\coordinate (v17) at (211,175);
\coordinate (v18) at (171,198);

\draw (v0) -- (v1);
\draw (v0) -- (v7);
\draw (v0) -- (v8);
\draw (v0) -- (v18);

\draw (v1) -- (v2);
\draw (v1) -- (v9);
\draw (v1) -- (v17);

\draw (v2) -- (v3);
\draw (v2) -- (v10);
\draw (v2) -- (v16);

\draw (v3) -- (v4);
\draw (v3) -- (v7);
\draw (v3) -- (v11);

\draw (v4) -- (v5);
\draw (v4) -- (v12);
\draw (v4) -- (v18);

\draw (v5) -- (v6);
\draw (v5) -- (v13);
\draw (v5) -- (v17);

\draw (v6) -- (v7);
\draw (v6) -- (v14);
\draw (v6) -- (v16);

\draw (v7) -- (v15);

\draw (v8) -- (v11);
\draw (v8) -- (v13);
\draw (v8) -- (v16);

\draw (v9) -- (v12);
\draw (v9) -- (v13);
\draw (v9) -- (v14);

\draw (v10) -- (v13);
\draw (v10) -- (v15);
\draw (v10) -- (v18);

\draw (v11) -- (v14);
\draw (v11) -- (v17);

\draw (v12) -- (v15);
\draw (v12) -- (v16);

\draw (v14) -- (v18);
\draw (v15) -- (v17);

\foreach \i in {0,...,18}{
    \pgfmathtruncatemacro{\labelnum}{\i+1}
    \node at (v\i) {$\labelnum$};
}

\end{tikzpicture}
\caption{The Robertson graph}
\label{P17.fig1}
\end{figure}
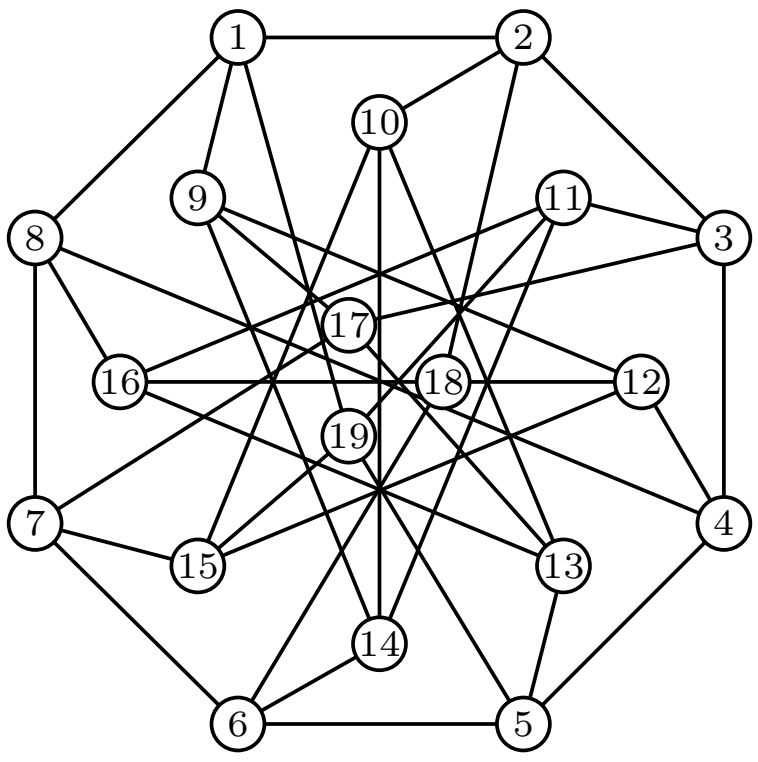
\begin{theorem}
Let $\mathcal R$ be the Robertson graph. Then $\mu(\mathcal R)=9$.
\end{theorem}
\begin{proof}
Since $\mathcal R$ is a $(4,5,+2)$-graph \cite{Rob_1964},
Theorem~\ref{P17.th3} gives $\mu(\mathcal R)\leq 9$. Consider the vertex set
$S=\{1,3,5,7,9,11,13,15,18\}$ in Figure~\ref{P17.fig1}.
Then $\mathcal R[S]\cong 4K_2\cup K_1$. Since $\Delta(\mathcal R[S])=1$, every pair
of vertices of $S$ at distance at most $2$
is $S$-visible.

The only pairs of vertices of $S$ at
distance greater than $2$ are
$\{1,13\}$, $\{3,15\}$, $\{5,9\}$, and $\{7,11\}$. These pairs are joined by the respective
shortest paths $(1,2,10,13), (3,2,10,15), (5,4,12,9), (7,6,14,11)$ whose internal vertices lie outside $S$. Hence $S$ is a mutual-visibility set of
$\mathcal R$, and result follows.
\end{proof}

\begin{table}[h]
\centering
\caption{Mutual-visibility numbers and bounds for
cages of girth $5$.}
\label{P17.tab1}
\begin{tabular}{lcccc}
\hline
Graph & $d$ & $\delta$ & $\mu(G)$ or bound & Status \\
\hline
Petersen graph
    & 3 & 0 & $6$ & Settled~\cite{Cice_2023} \\
Robertson graph
    & 4 & 2 & $9$ & This paper \\
$(5,5)$-cages
    & 5 & 4 & $\leq 15$ & Open \\
$(6,5)$-cage
    & 6 & 3 & $\leq 18$ & Open \\
Hoffman-Singleton graph
    & 7 & 0 & $20$ & Settled~\cite{TonShi_dif2026} \\
\hline
\end{tabular}
\end{table}
\section*{Acknowledgements}
The author is grateful to Dr. Shikhi M. for her encouragement and
motivation throughout this work, and for carefully reading the
manuscript and providing valuable comments and suggestions.
\bibliographystyle{plainurl}
\bibliography{cas-refs}
\end{document}